\documentclass[12pt,a4paper]{article}

\usepackage{amsmath,amssymb,epsf}
\usepackage[T1, T2A]{fontenc}
\usepackage[utf8]{inputenc}
\usepackage[russian, english]{babel}

\def\ot{\otimes}

\def\cA{{\cal A}}

\def\cG{{\cal G}}
\def\cI{{\cal I}}
\def\cH{{\cal H}}

\def\cS{{\cal S}}
\def\cW{{\cal W}}

\def\cU{{\cal U}}

\def\Aut{{\rm Aut}}

\def\N{{\mathbb N}}

\def\C{{\mathbb C}}

\def\d{{\partial}}

\usepackage{theorem}

\newtheorem{theorem}{Theorem}[section]
\newtheorem{proposition}{Proposition}[section]

{\theorembodyfont{\rmfamily}
\newtheorem{definition}[proposition]{Definition}
\newtheorem{example}[proposition]{Example}
\newtheorem{remark}[proposition]{Remark}

}
\newcommand{\rb}{\raisebox}
\newcommand{\ig}{\includegraphics}
\newcommand\risS[6]{\rb{#1pt}[#5pt][#6pt]{\begin{picture}(#4,15)(0,0)
  \put(0,0){\ig[width=#4pt]{#2.eps}} #3
     \end{picture}}}

\usepackage{graphicx}

\title{Polynomial graph invariants and the KP hierarchy}

\author{
Sergei Chmutov\thanks{Ohio State University, chmutov@math.osu.edu},
Maxim Kazarian\thanks{Steklov Mathematical Institute RAS,
National Research University Higher School of Economics,
Skolkovo Institute of Science and Technology, kazarian@mccme.ru},
Sergey Lando\thanks{National Research University Higher School of Economics,
Skolkovo Institute of Science and Technology, lando@hse.ru}}

\date{}

\begin{document}
\maketitle

\begin{abstract}
We prove that the generating function for the symmetric chromatic polynomial of all connected graphs
satisfies (after appropriate scaling change of variables) the Kadomtsev--Petviashvili
integrable hierarchy of mathematical physics. Moreover, we describe a large family of
polynomial graph invariants giving the same solution of the KP.
In particular, we introduce the Abel polynomial for graphs and show this for its generating function.
The key point here is a Hopf algebra structure on the space spanned by graphs and the behavior of the invariants on its primitive space.
\end{abstract}

\section{Introduction}
{\it\bfseries The Kadomtsev--Petviashvili {\rm(}KP{\rm)} hierarchy} is an infinite system of nonlinear partial differential equations for a function $F(p_1,p_2,\dots)$ of infinitely many variables, see Sec.~\ref{sIH} for details. The equations are indexed by partitions
of integers $n$, $n\ge4$, into two parts none of which is 1.
The first two equations, those corresponding to partitions of 4 and 5, are
$$\frac{\d^2F}{\d p_2^2} = \frac{\d^2F}{\d p_1\d p_3}-
      \frac12\Bigl(\frac{\d^2F}{\d p_1^2}\Bigr)^2 - \frac1{12}\frac{\d^4F}{\d p_1^4}$$
$$\frac{\d^2F}{\d p_2\d p_3} = \frac{\d^2F}{\d p_1\d p_4}-
      \frac{\d^2F}{\d p_1^2}\cdot\frac{\d^2F}{\d p_1\d p_2}
			- \frac16\frac{\d^4F}{\d p_1^3\d p_2}$$
The left hand side of the equations correspond to partitions into two parts none of which is 1,
while the terms on the right hand sides correspond to partitions of the same number $n$ involving 1.
For $n=6,$ there are two equations, which correspond to the partitions $2+4=6$ and
$3+3=6$, and so on.

It is well known now~\cite{O00} that the generating function for simple connected Hurwitz numbers
is a solution to the Kadomtsev--Petviashvili hierarchy.
Simple connected Hurwitz numbers enumerate certain connected embedded graphs, and
similar statements are known to be true for generating
functions enumerating  other kinds of embedded graphs as well.
In the present paper, we derive a similar statement for the generating functions for
Stanley's symmetric chromatic polynomial of graphs,
the Abel polynomials enumerating rooted forests in graphs,
and for a more general family of polynomial graph invariants.
As far as we know, this is the first appearance of KP solutions coming from graph polynomials.

\bigskip
{\it\bfseries The symmetric chromatic polynomial.}
In 1995, Richard Stanley~\cite{St95} introduced the {\it symmetric chromatic polynomial}
$X_G(x_1,x_2,\dots)$ of a finite graph $G$ with the vertex set $V(G)$
as the sum over all proper colorings $\kappa:V(G)\to\N$ that color the endpoints of each edge into distinct colors:
$$X_G(x_1,x_2,...) = \sum_{\kappa} \prod_{v\in V(G)} x_{\kappa(v)}.
$$
Being a symmetric function in the variables $x_1,x_2,\dots$,
the function~$X_G$ can be expressed through a basis of symmetric functions.

In fact, a polynomial $W_G(q_1,q_2,\dots)$, equivalent to the expression of $X_G(p_1,p_2,...)$ in terms of the power-sum functions
$p_m := \sum_{i=1}^{\infty} x_{i}^{m}$, was introduced in 1994 in \cite{CDL94}, motivated by some problems in knot theory.
The equivalence was noted in \cite[Theorem 6.1]{NW99}: $$(-1)^{|V(G)|}W_G(q_j=-p_j)=X_G(p_1,p_2,...).$$

Following \cite{CDL94}, we call $W_G(q_1,q_2,\dots)$ the {\it weighted chromatic polynomial}.
\label{pW_G}
There are two advantages of using the \cite{CDL94} approach as opposed to the \cite{St95} one.
Firstly, $W_G(q_1,q_2,\dots)$ satisfies a contraction-deletion relation, which is understood properly
in terms of graphs with weighted vertices. Secondly, $W_G(q_1,q_2,\dots)$ encodes the graph $G$ as an element of an appropriate Hopf algebra in terms of its primitive elements. We give the precise definitions and review the relevant results in Sec.\ref{sec:HA-graphs}.

The main result of this paper is that after an appropriate rescaling of the variables,  the function
$$\cW(q_1,q_2,\dots):=
\sum_{G\text{ connected}}\frac{W_G(q_1,q_2,\dots)}{|\Aut(G)|}
$$
becomes a solution to the KP hierarchy. Here and below $|\Aut(G)|$ denotes the order of the automorphism group of the graph $G$.
Note that we do not include the empty graph in the set of connected graphs.

\bigskip
{\it\bfseries The Abel polynomial.}
The classical sequence of Abel polynomials is defined as $x(x-an)^{n-1}$.
We consider its specialization with the parameter $a=-1$, $A_n(x):=x(x+n)^{n-1}$.
This is a sequence of polynomials of degree $n=\deg A_n$, $n=0,1,2,\dots$,
with leading coefficient~$1$ and zero free term (for $n>0$).
It possesses the following remarkable properties:
\begin{itemize}
\item $\frac{d}{dx} e^{-d/dx} (A_n(x)) = nA_{n-1}(x)$, see \cite{RST};
\item it is of {\it binomial type} \cite{Ab826}, $A_n(x+y)=\sum_{k=0}^{n} {n\choose k}A_k(x)A_{n-k}(y)$;
\item its coefficient of $x^k$ is equal to the number of rooted spanning forests
 consisting of $k$ trees (each of which has a root) in the complete graph $K_n$. For instance, the coefficient of $x$ equals $n^{n-1}$, the number of rooted spanning trees of $K_n$, according to the famous Cayley formula.
\end{itemize}

The last property allows us to define the Abel polynomial of an arbitrary graph
and make it into a polynomial of several variables $q_1,q_2,\dots$. Namely, with
each spanning forest of a graph $G$, we associate a monomial equal to the product of
$(i\cdot q_i)$ over the trees of the forest, where $i$ stands for the number of vertices of the tree, and the multiplication by $i$ corresponds to the different choices of the root of the tree. The {\it Abel polynomial} $A_G(q_1,q_2,\dots)$ of $G$
 is equal to the sum of these monomials over all (non-rooted) spanning forests of $G$.

 In particular, for a graph $G$ with $n$ vertices, the coefficient of $q_n$ in $A_G$ is equal to the number of rooted spanning trees of $G$. For a complete graph $K_n$, after the substitution
$q_i=x$ for all $i$, the polynomial becomes the Abel polynomial $A_n(x)$, $A_{K_n}(x,x,\dots)=A_n(x)$.
This polynomial is related to the spectrum of the combinatorial Laplace operator via the matrix-forest theorem \cite{CS06,Kn13}.

%

An analog of the second property of the Abel polynomial relates it to the Hopf algebra structure on graphs.
We discuss it in detail in Sec.\ref{sec:HA-graphs}. As a consequence, we get one more
result: after an appropriate rescaling of the variables,  the function
$$\cA(q_1,q_2,\dots):=\sum_{G\text{ connected}}\frac{A_G(q_1,q_2,\dots)}{|\Aut(G)|}
$$
becomes a solution to the KP hierarchy.

\bigskip
The work on this paper started during the second author's visit to the Ohio State University,
to which he expresses his gratitude.
The second and the third authors appreciate the support of RSF grant, project
16-11-10316 dated 11.05.2016.

\section{Main results}\label{sI}

{\it\bfseries The weighted chromatic polynomial.}\\
Necessary definitions will be given in Sec.~\ref{sec:HA-graphs}.
Besides the function $\cW(q_1,q_2,\dots)$, we introduce the function
$$
\cW^\circ(q_1,q_2,\dots):=\exp(\cW)=\sum_{G}\frac{W_G(q_1,q_2,\dots)}{|\Aut(G)|}
$$
represented by the sum over all (not necessarily connected, possibly empty) simple graphs.

For rescaling the variables we need the sequence \cite{Sloane}
$$
c_1=1;\qquad
c_n=(-1)^{n+1}\left[1+ \sum_{k=1}^{n-1} (-1)^k 2^{k(n-k)} {n-1\choose k-1} c_k\right]\ .
$$
In particular,
$c_1=1,\quad c_2=1,\quad c_3=5,\quad c_4=79,\quad c_5=3377,\dots
$.

\begin{theorem}\label{twc}
After the substitution $\displaystyle q_n=\frac{2^{n(n-1)/2}(n-1)!}{c_n}\cdot p_n$,
the generating function $\cW$ becomes a solution to the KP hierarchy of partial differential equations and $\cW^\circ$ becomes a $\tau$-function of the KP hierarchy.
\end{theorem}

{\it\bfseries The Abel polynomial.}\\
Similarly, we introduce the function
$$\cA^\circ(q_1,q_2,\dots):=\exp(\cA)=\sum_G\frac{A_G(q_1,q_2,\dots)}{|\Aut(G)|},
$$
and the constants $a_n:=2^{(n-1)(n-2)/2}n^{n-1}$ for $n=1,2,\dots$, so that $a_1=1$, $a_2=2$, $a_3=18$, $a_4=512$.

\begin{theorem}\label{ta}
After  the substitution $\displaystyle q_n=\frac{2^{n(n-1)/2}(n-1)!}{a_n}\cdot p_n$, the generating function $\cA$ becomes the solution to the KP hierarchy from Theorem \ref{twc}. So $\cA^\circ$ becomes the $\tau$-function of the KP hierarchy from the same theorem.
\end{theorem}

Theorems \ref{twc} and \ref{ta} are corollaries of the following more general statement
concerning umbral graph polynomial invariants as defined in Sec.~\ref{sec:HA-graphs}.
Both the weighted chromatic polynomial and the Abel polynomial are examples of umbral graph invariants.

\begin{theorem}\label{tm}
Let~$I$ be an umbral graph polynomial invariant with values in the polynomial ring in infinitely many variables $q_1,q_2,\dots$.
Define two generating functions by
\begin{eqnarray*}
\cI^\circ(q_1,q_2,\dots)&=&\sum_{G}\frac{I_G(q_1,q_2,\dots)}{|\Aut(G)|}\\
\cI(q_1,q_2,\dots)&=&\sum_{G\text{\rm\ connected}}\frac{I_G(q_1,q_2,\dots)}{|\Aut(G)|}
\end{eqnarray*}
where the first summation is carried over all graphs, while the second one over
all connected graphs and $|\Aut(G)|$ denotes the order of the automorphism group
of the graph~$G$.

Then after an appropriate rescaling of the variables~$\displaystyle q_n=\frac{2^{n(n-1)/2}(n-1)!}{i_n}\cdot p_n$, where~$i_n$, $n=1,2,\dots$, are some constants,
the generating function~$\cI$ becomes a solution to the KP hierarchy of partial differential equations.
After the same rescaling of the variables the generating function~$\cI^\circ$ becomes a
$\tau$-function of the KP hierarchy.

The rescaling constants $i_n$ are given by the sum
$$
i_n=n!\sum_{G\text{\rm\ connected}} \frac{[q_n]I_G(q_1,q_2,\dots)}{|\Aut(G)|}
$$
of the coefficients of~$q_n$.
\end{theorem}

We would like to stress that, after the rescaling in the theorem, any umbral graph invariant $I$ leads
to one and the same generating function.

In the last formula the summation can well be carried over all graphs rather than only connected ones,
since the coefficient of~$q_n$ in the value $I_G(q_1,q_2,\dots)$
on any disconnected graph~$G$ with~$n$ vertices is~$0$.
Note also that the value of an umbral graph invariant~$I_G$ on any graph with~$n$ vertices is
a weighted homogeneous polynomial of weight~$n$, with weights of the variable~$q_i$
taken to be~$i$, for $i=1,2,\dots$.

The previous theorem can be made more precise.


\begin{theorem}\label{tef}
After the rescaling of the variables described in Theorem~{\rm\ref{tm}}, the generating function~$\cI^\circ$
becomes the following linear combination of one-part Schur polynomials:
$$
\cS(p_1,p_2,\dots)=
1+2^0s_1(p_1)+2^1s_2(p_1,p_2)+\dots+2^{n(n-1)/2}s_n(p_1,p_2,\dots,p_n)+\dots.
$$
\end{theorem}

For the discussion of Schur polynomials, see Sec.~\ref{sS} below.

\section{The Hopf algebra of graphs and umbral graph invariants}\label{sec:HA-graphs}

In this section we describe the Hopf algebra structure on the space spanned by graphs.
The study of Hopf algebras in combinatorics was initiated in \cite{JR79}. It was then developed in \cite{Sch94, Sch95}.
Later it led to the concept of combinatorial Hopf algebra (see \cite{ABS06} and \cite{BM16}) suitable for quasisymmetric functions.
Independently the Hopf algebra of graphs was introduced in \cite{CDL94}.

Let~$\cG$ denote the vector space over~$\C$ spanned by simple graphs. This vector space is graded,
$$
\cG=\cG_0\oplus\cG_1\oplus\cG_2\oplus\dots,
$$
where~$\cG_n$ denotes the finite dimensional vector space spanned by all graphs with~$n$ vertices.
The product of two graphs is just their disjoint union, while the value of the coproduct $\Delta:\cG\to\cG\ot\cG$
on a graph $G$ is obtained as the sum of all terms $G(V_1)\ot G(V_2)$ where $G(V_1)$ and $G(V_2)$ are two induced subgraphs of $G$ whose vertex sets $V_1$ and $V_2$ constitute
an ordered partition of the vertex set of $G$, $V(G)=V_1\sqcup V_2$.
Clearly, both the product and the coproduct respect the grading.
These operations make the vector space spanned by graphs into a {\it graded bialgebra}.
This bialgebra is commutative, cocommutative, with a unit given by the empty graph.
According to the classical Milnor--Moore theorem it is isomorphic to the algebra of polynomials over its primitive vector space,
which is the space of elements $x$ on which the comultiplication $\Delta$ acts as $\Delta(x)=x\ot1+1\ot x$.
Any such bialgebra can be thought of as a Hopf algebra with the antipode acting by multiplication by $-1$ on the primitive space and counit taking value $1$ on the empty graph and $0$ on all other graphs.

In the case of $\cG$, the reference to the Milnor--Moore theorem can be eliminated by the following observation.
Each homogeneous subspace~$\cG_n\subset\cG$ admits a representation as a direct sum $\cG_n=P(\cG_n)\oplus D(\cG_n)$
of the subspace $P(\cG_n)$ of primitive elements and the subspace $D(\cG_n)$ of decomposable elements, which is spanned
by products of elements having grading smaller than~$n$.

Now, the projection $\pi:\cG_n\to P(\cG_n)$ to the primitive space along the space of decomposable
elements is given by
\begin{equation}\label{epp}
\pi(G)=\sum_{V(G)=\!\bigsqcup\limits_{\beta\in B}\!\!\! V_\beta}(-1)^{|B|-1}(|B|-1)!\prod_{\beta\in B} G(V_\beta),
\end{equation}
where the summation is carried over all unordered partitions $B$ of the set~$V(G)$ of vertices of~$G$ into disjoint
unions of nonempty subsets, $|B|$ is the number of parts, and
$G(V_\beta)$ denotes the subgraph of~$G$ induced by the subset $V_\beta\subset V(G)$ (see~\cite{L97}, \cite{L00}). Note that the product $\prod_{\beta\in B} G(V_\beta)$ is represented by the graph obtained from $G$ by removing all edges connecting vertices belonging to different parts of the partition.

W.~Schmitt's interpretation~\cite{Sch94} of the projection to the subspace of primitive elements
as the logarithm of the identity mapping allows one to invert Eq.~(\ref{epp})
and express any graph as a polynomial in the primitive elements of the form $\pi(G)$:
\begin{equation}\label{epd}
G=\sum_{V(G)=\!\bigsqcup\limits_{\beta\in B}\!\!\! V_\beta}\;\prod_{\beta\in B} \pi(G(V_\beta)).
\end{equation}

The algebra of polynomials $\C[q_1,q_2,\dots]$ is another commutative and cocommutative Hopf algebra.
The multiplicaton acts as the usual multiplication of polynomials and the
comultiplication $\Delta$ acts on variables $q_n$ as on primitive elements,
$\Delta(q_n)=q_n\ot1+1\ot q_n$. It becomes graded if we set the degree of $q_n$ to be $n$.
Thus it has precisely one-dimensional primitive space spanned by $q_n$ in each degree $n$.

\begin{definition}
An {\it umbral graph polynomial invariant} is an arbitrary homomorphism $F:\cG\to\C[q_1,q_2,\dots]$
of graded Hopf algebras.
\end{definition}

The term `umbral' is justified by the interpretation of the classical umbral calculus
as the study of the Hopf algebra of polynomials given in~\cite{RR78}.
Below, the adjective `polynomial' will be often omitted since we consider no
other kinds of graph invariants.
Umbral invariants participate in Theorem~\ref{tm}.

Equivalently, umbral invariants can be defined as multiplicative invariants taking primitive elements
of~$\cG$ to linear polynomials (the latter are the primitive elements in the
Hopf algebra of polynomials).
The condition of respecting comultiplication by $F$ can be reformulated in terms of {\em the binomial property\/}, saying that for any  graph~$G$,
$$
F_G(x_1+y_1,x_2+y_2,x_3+y_3,\dots)=\sum_{V(G)=V_1\sqcup V_2}
F_{G(V_1)}(x_1,x_2,\dots)F_{G(V_2)}(y_1,y_2,\dots),
$$
where the summation on the right runs over all partitions of the set of vertices~$V(G)$
into an ordered disjoint union of two subsets.
So any umbral graph invariant is uniquely determined by a set of constants $b_G$, one for each graph $G$, such that $F_{\pi(G)}=b_G q_{|V(G)|}$.

\bigskip
{\it\bfseries The weighted chromatic polynomial.}\\
The weighted chromatic polynomial is defined in terms of one more graded Hopf algebra, introduced in \cite{CDL94}, namely,
the Hopf algebra of {\em weighted graphs} $\cH$ (it was denoted $\cW$ in \cite{CDL94}, but in the present paper we reserve
the notation $\cW$ for the generating function) together with a graded homomorphism $\cG\to\cH$. Also it was proved there that
$\cH$ is isomorphic to $\C[q_1,q_2,\dots]$. Hence we get an umbral graph invariant which makes Theorem~\ref{twc} about the generating function for weighted chromatic polynomials a corollary of the general Theorem~\ref{tm}.

The algebra $\cH$ is the quotient of the space spanned by {\it weighted graphs\/} modulo the
{\it weighted chromatic relation}.

\begin{definition} A {\em weighted graph} is a graph $G$ without loops and multiple edges given together with
a {\em weight} $w:V(G)\to\N$ that assigns a positive integer to each vertex of the graph.
\end{definition}
Ordinary simple graphs can be treated as
weighted graphs with the weights of all vertices equal to 1.

\begin{definition} The {\em weighted chromatic relation\/} is
the set of relations $\displaystyle G = G'_e + G''_e,$ one for each graph~$G$ and an edge~$e$ in it,
where the graph $G'_e$ is obtained from $G$ by removing the edge $e$; and $G''_e$ is obtained from $G$ by a contraction of $e$ such that
if a multiple edge arises, it is reduced to a single edge and the weight $w(v)$ of the new vertex $v$ is set up to be equal to the sum of the weights of the two ends of the edge $e$; weights of the other vertices remaining unchanged.
\end{definition}

The quotient space of the space spanned by all weighted graphs modulo these relations inherits
the multiplication (disjoint union) and the comultiplication (splitting the vertex set into two subsets) and therefore becomes a Hopf algebra. This is $\cH$. It is graded by the total weight of a graph, $\displaystyle w(G)=\sum_{v\in V(G)}w(v)$, so that
$$
\cH=\cH_0\oplus\cH_1\oplus\cH_2\oplus\dots,
$$
where the vector space~$\cH_n$ is spanned by all weighted graphs of total weight~$n$ modulo weighted chromatic relations.

It was proved in \cite{CDL94} that the Hopf algebra $\cH$ has a one-dimensional primitive space in each grading and thus is isomorphic to $\C[q_1,q_2,\dots]$. Iterating the weighted chromatic relation we can reduce any
graph to a linear combination of (disconnected) graphs without edges. In degree $n$, the primitive space
$P(\cH_n)$ is spanned by the graph which is a single vertex of weight $n$. If we send it to
$q_n$, then the image of an ordinary graph $G$ (considered as a weighted graph with weights of all vertices equal to 1) in $\cH$ can be represented by a polynomial in the variables $q_n$. This is exactly the {\it weighted chromatic polynomial} $W_G(q_1,q_2,\dots)$ from the Introduction (p.~\pageref{pW_G}).

A closed formula for $W_G(q_1,q_2,\dots)$ is found in \cite[Corollary 4.4]{NW99},
\begin{equation}\label{esscp}
W_G(q_1,q_2,\dots)=\sum_{E'\subseteq E(G)} q_{v_1}\dots q_{v_k} (-1)^{|E'|-|V(G)|+k(E')}\ ,
\end{equation}
where the sum runs over all subsets $E'$ of the edge set $E(G)$, which we may consider as
spanning subgraphs of $G$, $k(E')$ is the number of connected components of this spanning subgraph,
$v_1,\dots, v_k$ denoting the numbers of vertices in the connected components.

\begin{example}\label{ex:wch}\

\underline{$\mathbf{n=1}.$} There is only one graph with one vertex, $G=\bullet$. It is primitive of degree 1 as a weighted graph,
whence $W_\bullet(q_1,q_2,\dots)=q_1$.

\bigskip
\underline{$\mathbf{n=2}.$} There is only one connected (simple) graph with two vertices,
$G=\risS{1}{gr_A2}{}{20}{0}{0}$. The order of its automorphism group is~$2$. According to the
weighted chromatic relation, it can be represented (as a weighted graph) as
$G=\risS{1}{gr_A2}{}{20}{0}{0}=(\bullet\ \bullet) +
\begin{array}[t]{c}\bullet\\[-8pt] \scriptstyle 2\end{array}$,
where the weight of a vertex is written underneath, weights $1$ being omitted.
The first graph is disconnected, whence taken to $q_1^2$ in $\C[q_1,q_2,\dots]$. The second graph is sent to $q_2$. Therefore,
$W_{\risS{-3}{gr_A2}{}{20}{0}{0}}(q_1,q_2,\dots)=q_1^2+q_2$.

\bigskip
\underline{$\mathbf{n=3}.$} There are two connected simple graphs with three vertices, namely,
$\risS{0}{gr_A3}{}{35}{0}{0}$ and $\risS{-5}{gr_tr}{}{20}{15}{0}$. They have automorphism groups of orders 2 and 6, respectively. As weighted graphs, modulo the weighted chromatic relation, they can be reduced as follows:
$$\begin{array}{rcl}
\risS{0}{gr_A3}{}{35}{0}{0} &=& (\risS{0}{gr_A3-e}{}{35}{0}{0}) +
			\begin{array}[t]{l}\risS{1}{gr_A2}{}{20}{0}{0}\\[-8pt] \scriptstyle 2\end{array} =
	    (\risS{0}{gr_A3-ee}{}{35}{0}{0}) +
			2(\begin{array}[t]{r}\bullet\ \bullet\\[-8pt] \scriptstyle 2\end{array})+
      \begin{array}[t]{c}\bullet\\[-8pt] \scriptstyle 3\end{array}\\[8pt]
\risS{-5}{gr_tr}{}{20}{15}{0} &=& \risS{-5}{gr_tr-e}{}{20}{15}{0}+
			\begin{array}[t]{l}\risS{1}{gr_A2}{}{20}{0}{0}\\[-8pt] \scriptstyle 2\end{array} =
		\Bigl(\risS{-5}{gr_tr-eee}{}{20}{15}{0}\Bigr) +
			3(\begin{array}[t]{r}\bullet\ \bullet\\[-8pt] \scriptstyle 2\end{array})+
      2\begin{array}[t]{c}\bullet\\[-8pt] \scriptstyle 3\end{array}
\end{array}
$$
The corresponding weighted chromatic polynomials are therefore
$$W_{\risS{-3}{gr_A3}{}{35}{0}{0}}(q_1,q_2,\dots)=q_1^3+2q_1q_2+q_3, \qquad
W_{\!\!\!\risS{-12}{gr_tr}{}{20}{0}{0}\!\!\!}(q_1,q_2,\dots)=q_1^3+3q_1q_2+2q_3.
$$

\underline{$\mathbf{n=4}.$} There are six connected simple graphs with four vertices. We arrange the results in the following table.
$$\begin{array}{|c|c|c|}\hline
G & \makebox(0,15){\ } W_G(q_1,q_2,\dots) & |\Aut(G)| \\[5pt] \hline\hline
\risS{0}{gr_A4}{}{55}{13}{10} & q_1^4+3q_1^2q_2+q_2^2+2q_1q_3+q_4 & 2 \\ \hline
\risS{-9}{gr_D4}{}{35}{16}{14} & q_1^4+3q_1^2q_2+3q_1q_3+q_4 & 6 \\ \hline
\risS{-9}{gr_D4+e}{}{35}{16}{14} & q_1^4+4q_1^2q_2+q_2^2+4q_1q_3+2q_4 & 2\\ \hline
\risS{-8}{gr_sq}{}{22}{20}{20} & q_1^4+4q_1^2q_2+2q_2^2+4q_1q_3+3q_4 & 8 \\ \hline
\risS{-8}{gr_sq+d}{}{22}{20}{20} & q_1^4+5q_1^2q_2+2q_2^2+6q_1q_3+4q_4 & 4 \\ \hline
\risS{-8}{gr_K4}{}{22}{20}{20} & q_1^4+6q_1^2q_2+3q_2^2+8q_1q_3+6q_4 & 24 \\ \hline
\end{array}
$$
Consequently, the generating function
$$\cW(q_1,q_2,\dots)=
\sum_{G{\text{ connected}}\atop{\text{\quad\ non-empty}}}\frac{W_G(q_1,q_2,\dots)}{|\Aut(G)|}\\
$$
starts with
\begin{eqnarray*}
\cW(q_1,q_2,\dots)&=&
\frac1{1!}q_1+\frac1{2!}\left(q_1^2+q_2\right)+\frac1{3!}\left(4q_1^3+9q_1q_2+5q_3\right)\\
&&+\frac1{4!}
\left(38q_1^4+144q_1^2q_2+45q_2^2+140q_1q_3+79q_4\right)+\dots,
\end{eqnarray*}

The beginning of its exponent looks like
\begin{eqnarray*}
\cW^\circ(q_1,q_2,\dots)&=&\exp(\cW)=\sum_{G}\frac{W_G(q_1,q_2,\dots)}{|\Aut(G)|}\\
&=&1+\frac1{1!}q_1+\frac1{2!}\left(2q_1^2+q_2\right)+\frac1{3!}\left(8q_1^3+12q_1q_2+5q_3\right)\\
&&+\frac1{4!}
\left(64q_1^4+192q_1^2q_2+48q_2^2+160q_1q_3+79q_4\right)+\dots\\
\end{eqnarray*}

The rescaling of variables in Theorem \ref{twc} gives
$$q_1=p_1,\quad q_2=2p_2,\quad q_3=\frac{2^3\cdot2!}{5}p_3,\quad q_4=\frac{2^6\cdot3!}{79}p_4,\quad\dots
$$
After this substitution we get the function
\begin{eqnarray*}
F(p_1,p_2,\dots)&=&\cW(p_1,2p_2,\frac{16}{5}p_3,\frac{384}{79}p_4,\dots)\\ &=&
p_1+\frac12\left(p_1^2+2p_2\right)+\frac16\left(4p_1^3+18p_1p_2+16p_3\right)\\
&&+\frac1{24}
\left(38p_1^4+288p_1^2p_2+180p_2^2+448p_1p_3+384p_4\right)+\dots,
\end{eqnarray*}
which coincides with the logarithm of the function~$\cS$ from Theorem~\ref{tef}.

One may easily check that $\frac{\d^2F}{\d p_2^2} = 15 +\dots$,\qquad
while
$$\frac{\d^2F}{\d p_1\d p_3}-
      \frac12\Bigl(\frac{\d^2F}{\d p_1^2}\Bigr)^2 - \frac1{12}\frac{\d^4F}{\d p_1^4}
=\frac{56}3-\frac12-\frac{19}6+ \dots=15+ \dots$$
giving thus an evidence that we get a solution of the first KP equation.
\end{example}

\bigskip
{\it\bfseries The Abel polynomial.}\\

\begin{example}\label{ex:Abl} This example is similar to \ref{ex:wch}. We will follow it case by case.

\bigskip
\underline{$\mathbf{n=1}.$} For $G=\bullet$, we have $A_\bullet(q_1,q_2,\dots)=q_1$.

\bigskip
\underline{$\mathbf{n=2}.$} For $G=\risS{1}{gr_A2}{}{20}{0}{0}$, there are two spanning forests, $(\bullet\ \bullet)$ and the graph itself $\risS{1}{gr_A2}{}{20}{0}{0}$. Therefore,
$A_{\risS{-3}{gr_A2}{}{20}{0}{0}}(q_1,q_2,\dots)=q_1^2+2q_2$.

\bigskip
\underline{$\mathbf{n=3}.$} There are two connected simple graphs with three vertices,
$\risS{0}{gr_A3}{}{35}{0}{0}$ and $\risS{-5}{gr_tr}{}{20}{15}{12}$ ,with automorphism groups of order 2 and 6, respectively. For the first one, $\risS{0}{gr_A3}{}{35}{0}{0}$, there are four spanning forests
$$(\risS{0}{gr_A3-ee}{}{35}{0}{0}),\qquad (\risS{0}{gr_A3-e}{}{35}{0}{0}),\qquad
(\risS{0}{gr_A3-e-op}{}{35}{0}{0}),\qquad \risS{0}{gr_A3}{}{35}{0}{0}
$$
which contribute the following monomials to the Abel polynomial:
$$q_1^3,\qquad 2q_1q_2,\qquad 2q_1q_2,\qquad 3q_3.$$
Thus $A_{\risS{-3}{gr_A3}{}{35}{0}{0}}(q_1,q_2,\dots)=q_1^3+4q_1q_2+3q_3$.

For the triangle, $\risS{-5}{gr_tr}{}{20}{0}{0}$, there are seven spanning forests
$$\Bigl(\risS{-5}{gr_tr-eee}{}{20}{0}{0}\Bigr),\quad
\Bigl(\risS{-5}{gr_tr-ee}{}{20}{0}{0}\Bigr),\quad
\Bigl(\risS{-5}{gr_tr-ee-op}{}{20}{0}{0}\Bigr),\quad
\Bigl(\risS{-5}{gr_tr-ee-t}{}{20}{0}{0}\Bigr),\quad
\risS{-5}{gr_tr-e}{}{20}{0}{0},\quad
\risS{-5}{gr_tr-e-op}{}{20}{0}{0},\quad
\risS{-5}{gr_tr-e-t}{}{20}{0}{0}.
$$
Each of the last three spanning trees contributes $3q_3$. The previous three spanning forests contribute
$2q_1q_2$ each. Hence the Abel polynomial of the triangle is
$A_{\!\!\!\risS{-12}{gr_tr}{}{20}{0}{0}\!\!\!}(q_1,q_2,\dots)=q_1^3+6q_1q_2+9q_3$.

\bigskip
\underline{$\mathbf{n=4}.$} For the six connected simple graphs with four vertices, we collect the results in the following table.
$$\begin{array}{|c|c|c|}\hline
G & \makebox(0,15){\ } A_G(q_1,q_2,\dots) & |\Aut(G)| \\[5pt] \hline\hline
\risS{0}{gr_A4}{}{55}{13}{10} & q_1^4+6q_1^2q_2+4q_2^2+6q_1q_3+4q_4 & 2 \\ \hline
\risS{-9}{gr_D4}{}{35}{16}{14} & q_1^4+6q_1^2q_2+9q_1q_3+4q_4 & 6 \\ \hline
\risS{-9}{gr_D4+e}{}{35}{16}{14} & q_1^4+8q_1^2q_2+4q_2^2+15q_1q_3+12q_4 & 2\\ \hline
\risS{-8}{gr_sq}{}{22}{20}{20} & q_1^4+8q_1^2q_2+8q_2^2+12q_1q_3+16q_4 & 8 \\ \hline
\risS{-8}{gr_sq+d}{}{22}{20}{20} & q_1^4+10q_1^2q_2+8q_2^2+24q_1q_3+32q_4 & 4 \\ \hline
\risS{-8}{gr_K4}{}{22}{20}{20} & q_1^4+12q_1^2q_2+12q_2^2+36q_1q_3+64q_4 & 24 \\ \hline
\end{array}
$$
Thus the generating function
$$\cA(q_1,q_2,\dots)=
\sum_{G{\text{ connected}}\atop{\text{\quad\ non-empty}}}\frac{A_G(q_1,q_2,\dots)}{|\Aut(G)|}\\
$$
starts with
\begin{eqnarray*}
\cA(q_1,q_2,\dots)&=&
\frac1{1!}q_1+\frac1{2!}\left(q_1^2+2q_2\right)+\frac1{3!}\left(4q_1^3+18q_1q_2+18q_3\right)\\
&&+\frac1{4!}
\left(38q_1^4+288q_1^2q_2+180q_2^2+504q_1q_3+512q_4\right)+\dots\ .
\end{eqnarray*}

The rescaling of variables in Theorem \ref{ta} has the form
$$q_1=p_1,\quad q_2=p_2,\quad q_3=\frac89 p_3,\quad q_4=\frac34 p_4,\quad\dots
$$
After the substitution we come to the same series
\begin{eqnarray*}
F(p_1,p_2,\dots)&=&\cA(p_1,p_2,\frac89 p_3,\frac34 p_4,\dots)\\ &=&
p_1+\frac12\left(p_1^2+2p_2\right)+\frac16\left(4p_1^3+18p_1p_2+16p_3\right)\\
&&+\frac1{24}
\left(38p_1^4+288p_1^2p_2+180p_2^2+448p_1p_3+384p_4\right)+\dots,
\end{eqnarray*}
as before.
It is interesting to note that for an individual graph $G$, the rescalings of the
variables under consideration give different polynomials $W_G(p_1,2p_2,\frac{16}{5}p_3,\frac{384}{79}p_4,\dots)$ and
$A_G(p_1,p_2,\frac89 p_3,\frac34 p_4,\dots)$. For example, for the six graphs with four vertices from our tables, the coefficients of $p_4$ in $W_G$ will be $\frac{384}{79}, \frac{384}{79},
\frac{768}{79}, \frac{1152}{79}, \frac{1536}{79}, \frac{2304}{79}$, while in $A_G$ they will be
$3, 3, 9, 12, 24, 48$. However, after division by $|\Aut(G)|$ and summation to the generating functions they become equal:
\begin{eqnarray*}
\frac1{24}\cdot \frac{12\cdot384+4\cdot384+12\cdot768+3\cdot1152+6\cdot1536+2304}{79} =
\frac{384}{24} &=& 16 \\
\frac1{24}\cdot (12\cdot3+4\cdot3+12\cdot9+3\cdot12+6\cdot24+48)&=&16\ .
\end{eqnarray*}
\end{example}

\bigskip
\begin{theorem}
The Abel polynomial $A:\cG\to\C[q_1,q_2,\dots]$ is an umbral graph invariant.
\end{theorem}

It is clear that the Abel polynomial is homogeneous and multiplicative.
In order to prove that the mapping $A:\cG\to\C[q_1,q_2,\dots]$ respects comultiplication,
it suffices to show that the Abel polynomial possesses the binomial property, that is, that it satisfies
the relation
$$
A_G(x_1+y_1,x_2+y_2,
\dots)=\sum_{V(G)=V_1\sqcup V_2}A_{G(V_1)}(x_1,x_2,\dots)A_{G(V_2)}(y_1,y_2,\dots),
$$
where the summation on the right runs over all partitions of the set of vertices~$V(G)$
into an ordered disjoint union of two subsets.

The binomial property is satisfied because of the following argument. Consider a spanning forest
in the graph~$G$. The contribution of this spanning forest to the expression on the left is
$$
\prod_{i=1}^\infty (i(x_i+y_i))^{k_i},
$$
where~$k_i$ is the number of trees with~$i$ vertices in the forest, $k_i=0$ for~$i$ large enough.
These~$k_i$ trees can be distributed between the two subsets~$V'$, $V''$ of the partition
of the set of vertices~$V(G)$. There are ${k_i\choose\ell_i}$ ways to choose~$\ell_i$
trees for the first part thus leaving  $k_i-\ell_i$ trees for the second part, for $\ell_i=0,1,2,\dots,k_i$.
Hence the binomial property for the Abel polynomial follows from the usual binomial identity
$$
(x_i+y_i)^{k_i}=\sum {k_i\choose\ell_i} x_i^{\ell_i}y_i^{k_i-\ell_i}.
$$

\begin{remark}
Suppose we have a multiplicative polynomial graph invariant taking values in the ring
$\C[x]$ of polynomials in a single variable $x$. Suppose also that the value of this invariant on
any primitive element in the Hopf algebra of graphs is a linear polynomial,
that is, $x$ with some coefficient.
Then such a graph invariant admits a canonical {\em umbralization}.
Namely, we can associate to this graph invariant a homomorphism from the Hopf algebra
of graphs to the Hopf algebra $\C[q_1,q_2,\dots]$ of polynomials in infinitely many variables,
which is determined uniquely by the requirement that its value on the projection $\pi(G)$
of an arbitrary graph~$G$ to the subspace of primitive elements is $b_Gq_n$, where
$n=|V(G)|$ is the number of vertices of~$G$, and $b_G$ is the coefficient of~$x$ in the
value of the given invariant on~$G$.

For example, Stanley's symmetrized chromatic polynomial is the umbralization of the ordinary
chromatic polynomial. Similarly, the Abel polynomial we introduce is the umbralization of
the one-variable polynomial graph invariant whose coefficient of~$x^k$, for a given graph $G$,
is the number of rooted forests in $G$ consistsing of $k$ trees.
\end{remark}

\section{General information on the Kadomtsev--Petviashvili integrable hierarchy}\label{sIH}

In this section, we reproduce necessary information about the  Kadomtsev--Petviashvili
(KP) integrable hierarchy of partial differential equations.
A more detailed explanation, suitable for the first acquaintance, can be found, for example, in~\cite{KL15}.

\subsection{$\tau$-functions and solutions to the KP hierarchy}
The KP integrable hierarchy of partial differential equations is best described in terms of
its solutions and $\tau$-functions. The latter are the exponents of the former.
In 1983, Sato~\cite{S83} shown that the $\tau$-functions of the KP-hierarchy
are in a natural one-to-one correspondence with semi-infinite planes in
an infinite dimensional vector space of Laurent polynomials.

Let~$V$ be the vector space of Laurent polynomials in the variable~$z$.
The {\it semi-infinite Grassmannian\/} $G(\frac\infty2,V)$ (to be more precise,
an open dense cell in this Grassmannian)
consists of decomposable vectors in $P\Lambda^{\frac\infty2}V$,
i.e., of vectors of the form
$$
\beta_1(z)\wedge \beta_2(z)\wedge \beta_3(z)\wedge\dots,
$$
where $\beta_i$ are Laurent series in~$z$, and we have
$$
\beta_i(z)=z^{-i}+c_{i1}z^{-i+1}+c_{i2}z^{-i+2}+\dots.
$$

The {\it semi-infinite wedge product\/}
$\Lambda^{\frac\infty2}V$ of the vector space~$V$ is, by definition, the vector space spanned by the vectors
$$
v_\mu=z^{m_1}\wedge z^{m_2}\wedge z^{m_3}\wedge\dots,\qquad
m_1>m_2>m_3>\dots,\qquad m_i=\mu_i-i,
$$
where~$\mu$ is a partition, $\mu=(\mu_1,\mu_2,\mu_3,\dots)$,
$\mu_1\ge\mu_2\ge\mu_3\ge\dots$,
having all but finitely many parts equal to~$0$.
In particular, $m_i=-i$ for all~$i$ large enough.

The empty partition corresponds to the {\it vacuum vector\/}
$$
v_{\emptyset}=z^{-1}\wedge z^{-2}\wedge z^{-3}\wedge\dots.
$$
Similarly,
$$
v_{1^1}=z^{0}\wedge z^{-2}\wedge z^{-3}\wedge\dots,\quad
v_{2^1}=z^{1}\wedge z^{-2}\wedge z^{-3}\wedge\dots,
\quad v_{1^2}=z^{0}\wedge z^{-1}\wedge z^{-3}\wedge\dots,
$$
and so on.

Numbering of the basic vectors in the semi-infinite wedge product
$\Lambda^{{\frac\infty2}}V$ (the space of {\it fermions}) by partitions of
nonnegative integers establishes its natural vector space isomorphism with the space of {\it bosons},
that is, power series in infinitely many variables  $p_1,p_2,\dots$.
It is called the {\em boson-fermion isomorphism}.
This isomorphism takes a basic vector~$v_\mu$
to the Schur polynomial $s_\mu=s_\mu(p_1,p_2,\dots)$. The latter is a quaihomogeneous polynomial
of degree $|\mu|$ in the variables~$p_i$
(by definition, the degree of~$p_i$ is~$i$).
Schur polynomials form an additive basis of the vector space of power series in the variables~$p_i$.

\begin{definition}\rm
The {\it Hirota equations\/} are the Pl\"ucker equations of the embedding
of the semi-infinite Grassmannian into the projectivized semi-infinite wedge product
$P\Lambda^{\frac\infty2}V$. Solutions of the Hirota equations (that is,
semi-infinite planes) are called {\it $\tau$-functions} of the KP hierarchy.
\end{definition}

Let us represent a power series in~$p_1,p_2,\dots$ as a linear combination
$$
\sum_\nu c_\nu s_\nu(p_1,p_2,\dots)
$$
of Schur polynomials.
The Hirota equations are quadratic in the coefficients~$c_\nu$ of these linear combinations.

\begin{definition}\rm
The form the Hirota equations take for the logarithm of a $\tau$-function
are called the {\it Kadomtsev--Petviashvili hierarchy equations}.
\end{definition}

In other words, each $\tau$-function of the KP hierarchy can be obtained as the result of the following
procedure:
\begin{itemize}
\item take a semi-infinite plane $\beta_1(z)\wedge
\beta_2(z)\wedge\dots$ in~$V$;
 \item using the Laurent expansion of~$\beta_i$
 represent the corresponding point in the semi-infinite Grassmannian as a linear combination
 of the basic vectors~$v_\kappa$ and multiply by a constant in order to make~$1$ the coefficient
 of the vacuum-vector~$v_\emptyset$;
 \item replace each vector~$v_\kappa$ in the resulting linear combination
 with the corresponding Schur polynomial~$s_\kappa(p_1,p_2,\dots)$;
 this series in $p_1,p_2,\dots$ is the desired $\tau$-function.
\end{itemize}

Logarithms of $\tau$-functions are solutions to the KP hierarchy.

\subsection{Schur polynomials and constructions of combinatorial $\tau$-functions}\label{sS}

For~$\lambda$ being partitions of a given non negative integer~$n$, the Schur polynomials $s_\lambda$
form an additive basis in the vector space of symmetric
polynomials in infinitely many variables $x_1,x_2,\dots$.
In particular, the one-part Schur polynomial $s_{n^1}=s_n$ is defined as the complete homogeneous symmetric polynomial, which is the sum of all monomials of homogeneous degree $n$
in variables $x_1,x_2,\dots$.

After substituting for~$p_i$ the power sum functions
$$
p_i=x_1^i+x_2^i+\dots,\qquad i=1,2,\dots,
$$
the Schur polynomials turn into an additive basis in the vector space of quasi-homogeneous polynomials
in~$p_1,\dots,p_n$, with the weight of~$p_i$ taken to be~$i$, for $i=1,2,\dots$.

The expression of the one-part Schur polynomial $s_n$ in variables $p_1,p_2,\dots$ can be found from
the generating function
\begin{eqnarray*}
\sum_{n=0}^\infty s_{n}
&=&e^{\sum_{k=1}^\infty\frac{p_k}k}\\
&=&1+p_1+\frac1{2!}(p_1^2+p_2)+\frac1{3!}(p_1^3+3p_1p_2+2p_3)+\dots\ .
\end{eqnarray*}
For arbitrary partition~$\lambda=(\lambda_1,\lambda_2,\dots,\lambda_l)$, $\lambda_1\ge\lambda_2\ge\dots\ge\lambda_l$,
the Schur polynomial~$s_\lambda$ can be given by
the Jacobi--Trudi formula
$$
s_\lambda=\det\left(\begin{array}{ccccc}
s_{\lambda_1}&s_{\lambda_1+1}&s_{\lambda_1+2}&\dots&s_{\lambda_1+l-1}\\
s_{\lambda_2-1}&s_{\lambda_2}&s_{\lambda_1+1}&\dots&s_{\lambda_2+l-2}\\
\dots&\dots&\dots&\dots&\dots\\
s_{\lambda_l-l+1}&s_{\lambda_l-l+2}&s_{\lambda_l-l+3}&\dots&s_{\lambda_l}
\end{array}\right).
$$

The Schur polynomials constitute a convenient language for describing formal $\tau$-functions
as well as solutions to the KP hierarchy. In particular, the following statement is true.

\begin{proposition}\label{tlopc}
Any linear combination
$$
1+c_1s_1+c_2s_2+\dots
$$
of one-part Schur polynomials with constant coefficients $c_n$ is a $\tau$-function for the KP hierarchy.
\end{proposition}

Indeed, consider the semi-infinite plane $\beta_1(z)\wedge z^{-2}\wedge z^{-3}\wedge\dots$, where
$\beta_1(z)=z^{-1}+c_1 z^0+c_2z^1+c_3z^2+\dots$. Since
$$
\beta_1(z)\wedge z^{-2}\wedge z^{-3}\wedge\dots=v_\emptyset+c_1v_{1^1}+c_2v_{2^1}+\dots,
$$
this plane is taken, under the boson-fermion correspondence, to the power series in the Proposition.

\section{Proof of the main theorems}

In this section, we give the proofs of main theorems stated in Sec.~\ref{sI}.

\subsection{Proof of Theorem~\ref{tm}}

Define the {\it universal umbral graph invariant\/} $G\mapsto U_G(q_1,q_2,\dots)$ as the multiplicative graph invariant
taking value $b_Gq_{|V(G)|}$ on the primitive element $\pi(G)$, which is the projection of a
graph~$G$ to the subspace of primitive elements along the subspace of decomposable elements in~$\cG$,
for all connected graphs~$G$. The invariant $U_G$ takes values in the ring of polynomials in the
variables $q_1,q_2,\dots$ whose coefficients are polynomials in the variables~$b_H$, for all connected graphs~$H$.
(Of course, for a given graph~$G$, the coefficients
of the polynomial $U_G$ depend only on those indeterminates~$b_H$ such that $H$ is an induced subgraph of~$G$).
We call the invariant~$U$ universal because the projections $\pi(G)$ of connected graphs to
the subspace of primitive elements form a basis in the space of primitive elements of~$\cG$,
whence any umbral graph invariant with values in the ring $\C[q_1,q_2,\dots]$ can be induced
from~$U$ by choosing appropriate values of~$b_G$, for all connected graphs~$G$.

Introduce two generating functions
\begin{eqnarray*}
\cU^\circ(q_1,q_2,\dots)&=&\sum_{G}\frac{U_G(q_1,q_2,\dots)}{|\Aut(G)|}\\
\cU(q_1,q_2,\dots)&=&\sum_{G\text{ connected}}\frac{U_G(q_1,q_2,\dots)}{|\Aut(G)|}.
\end{eqnarray*}

As above, we have
$$
\cU=\exp(\cU^\circ).
$$

For $k=0,1,2,\dots,$ denote by
$$
u_k=\sum_{{G\text{ connected}\atop |V(G)|=k}}\frac{b_G}{|\Aut(G)|}
$$
the coefficient of the monomial~$q_k$ in either~$\cU$ or~$\cU^\circ$.
As well, denote by
\begin{eqnarray*}
\cU^\circ_k(q_1,q_2,\dots)&=&\sum_{G, |V(G)|=k}\frac{U_G(q_1,q_2,\dots)}{|\Aut(G)|}\\
\cU_k(q_1,q_2,\dots)&=&\sum_{{G\text{ connected}\atop |V(G)|=k}}\frac{U_G(q_1,q_2,\dots)}{|\Aut(G)|}
\end{eqnarray*}
the $k$~th homogeneous parts of the generating functions $\cU^\circ$ for all graphs
and $\cU$ for connected graphs, respectively, so that
\begin{eqnarray*}
\cU^\circ(q_1,q_2,\dots)&=&1+\cU^\circ_1(q_1,q_2,\dots)+\cU^\circ_2(q_1,q_2,\dots)+\dots\\
\cU(q_1,q_2,\dots)&=&\cU_1(q_1,q_2,\dots)+\cU_2(q_1,q_2,\dots)+\dots.
\end{eqnarray*}

Introduce notation
$$
Q_k(q_1,q_2,\dots)=s_k\left(\frac{u_1}{0!2^0}q_1,\frac{u_2}{1!2^1}q_2,\frac{u_3}{2!2^3}q_3,\dots\right),
$$
that is, $Q_k(q_1,q_2,\dots)$ is the result of the substitution $p_i=\frac{u_i}{(i-1)!\;2^{{i\choose 2}}}q_i$
into the one-part Schur function $s_k(p_1,p_2,\dots,p_k)$.

\begin{theorem}\label{th}
For each $k=0,1,2,\dots,$ we have
$$
\cU^\circ_k(q_1,q_2,\dots)=2^{{k\choose2}}Q_k(q_1,q_2,\dots,q_k),
$$
so that
$$
\cU^\circ(q_1,q_2,\dots)=Q_0+2^0Q_1+2^1Q_2+2^3Q_3+2^6Q_4+\dots.
$$
\end{theorem}

This means that after the specified rescaling of the variables,
$$
q_k=\frac{2^{{k\choose2}}(k-1)!}{u_k}p_k,\qquad
k=1,2,\dots,
$$
the generating function $\cU(q_1,q_2,\dots)$ becomes a linear combination of one-part
Schur functions, namely, the function~$\cS(p_1,p_2,\dots)$ from Theorem~\ref{tef}, whence, according to Proposition~\ref{tlopc},
a $\tau$-function for the KP hierarchy.

In order to prove the theorem, we will need the following statement.

Set~$b_G=0$ for disconnected graphs~$G$.

Now, Eq.~(\ref{epd}) shows that for any (either connected or disconnected) graph~$G$, we have
\begin{equation}\label{esp}
U_G(q_1,q_2,\dots)=\sum_{V(G)=\sqcup_\alpha V_\alpha}\prod_\alpha b_{G(V_\alpha)}q_{|V_\alpha|},
\end{equation}
where the summation is carried over all unordered partitions of the set~$V(G)$ of vertices of~$G$
into disjoint unions of subsets, and $G(V_\alpha)$ denotes the subgraph of~$G$ induced by the subset
$V_\alpha\subset V(G)$.

In order to prove Theorem~\ref{th}, let us replace summation over all graphs
with~$k$ vertices with summing over all spanning subgraphs of the complete graph~$K_k$
on~$k$ vertices. Obviously, for a graph~$G$ with~$k$ vertices we have
$$
\frac{U_G(q_1,q_2,\dots)}{|\Aut(G)|}=\frac1{k!}\sum_{{E\subset E(K_k)\atop  K_k|_E=G}}U_{K_k|_E}(q_1,q_2,\dots),
$$
where the summation is carried over all spanning subgraphs~$K_k|_E$ of~$K_k$ that are isomorphic to~$G$.
Hence,
$$
\cU^\circ_k(q_1,q_2,\dots)=\frac1{k!}\sum_{E\subset E(K_k)}U_{K_k|_E}(q_1,q_2,\dots),
$$
and we have replaced summation over isomorphism classes of graphs with~$k$ vertices with summation over subsets
in the set of edges of the complete graph with~$k$ vertices.

Making use of Eq.~(\ref{esp}) we obtain
$$
\cU^\circ_k(q_1,q_2,\dots)=\frac{1}{k!}\sum_{V(K_k)=\sqcup_\alpha V_\alpha}2^{{k\choose2}-
\sum_\alpha{|V_\alpha|\choose2}}\prod_\alpha \sum_{E\subset E(K_k(V_\alpha))}b_{K_k(V_\alpha)|_E}q_{|V_\alpha|}.
$$
Here ${k\choose2}-\sum_\alpha{|V_\alpha|\choose2}$ is the number of the edges of~$K_k$ connecting
vertices belonging to {\it distinct\/} subsets $V_\alpha$. Dividing both sides of the last equation by
$2^{{k\choose2}}$ we obtain
\begin{eqnarray*}
2^{-{k\choose2}}\cU^\circ_k(q_1,q_2,\dots)&=&
\frac{1}{k!}\sum_{V(K_k)=\sqcup_\alpha V_\alpha}\prod_\alpha2^{-
{|V_\alpha|\choose2}} |V_\alpha|! u_{|V_\alpha|}q_{|V_\alpha|}\\
&=&\frac1{k!}\sum_{V(K_k)=\sqcup_\alpha V_\alpha}\prod_\alpha (|V_\alpha|-1)! p_{|V_\alpha|}.
\end{eqnarray*}

The right-hand side of the last equation is independent of the variables
$b_G$ and coincides with the definition of the one-part Schur function~$s_k=s_k(p_1,p_2,\dots)$,
which proves the theorem.

\subsection{Specialization to weighted chromatic polynomial}

Theorem~\ref{twc} for the weighted chromatic polynomial~$W_G$ follows from the main theorem,
since the weighted chromatic polynomial is an example of umbral
graph invariant. For this invariant, the sequence
$$
1,1,5,79,3377,\dots
$$
of the values $c_1,c_2,c_3,\dots$ can be described as follows. The exponent of the exponential generating function
for the alternating sequence
$$
x-\frac1{2!2^1}x^2+\frac5{3!2^3}x^3-\frac{79}{4!2^6}x^4+\frac{3377}{5!2^{10}}x^5-\dots
$$
is
$$
1+x+\frac1{2!2^1}x^2+\frac1{3!2^3}x^3+\frac1{4!2^6}x^4+\frac1{5!2^{10}}x^5+\dots.
$$
In other words, we have the recursion
$$
(-1)^{n+1}c_n=1-\sum_{k=1}^{n-1}(-1)^{k+1}2^{k(n-k)}{n-1\choose k-1}c_k.
$$

In order to prove this statement, let us compute explicitly the coefficient $C_k$ of $x^k$ in the exponent
$$
\exp\left(\sum_{i=1}^\infty\frac{(-1)^{i-1}c_i}{i!2^{{i\choose 2}}}x^i\right).
$$
This exponent is the result of substitution $p_i=\frac{(-1)^{i-1}c_i}{(i-1)!\;2^{{i\choose 2}}}x^i$ into $\exp\sum_{i=1}^\infty(p_i/i)$, therefore, $C_kx^k$ is the value of the Schur function $s_k$ at the same substitution. By Theorem~\ref{th}, we get
$$
C_k=2^{-{k\choose 2}}\cW^\circ_k(q_1,q_2,\dots)\bigm|_{q_i=(-1)^{i-1}},
$$
where $\cW^\circ_k$ is the degree $k$ quasihomogeneous summand of~$\cW^\circ$. Now, following the proof of Theorem~\ref{th}, we get
$$
C_k=2^{-{k\choose 2}}\frac1{k!}\sum_{E\subset E(K_k)}W_{K_k|_E}(q_1,q_2,\dots)\bigm|_{q_i=(-1)^{i-1}},
$$
where the summation is carried over all spanning subgraphs~$K_k|_E$ of~$K_k$. Up to this point, similar computations 
can be applied to any umbral graph invariant. Now, for the particular case of the invariant~$W_G$, we get by~\eqref{esscp},
\begin{eqnarray*}
C_k&=&\frac{2^{-{k\choose 2}}}{k!}\sum_{E'\subset E\subset E(K_k)}(-1)^{|E'|}
=\frac{2^{-{k\choose 2}}}{k!}\sum_{r=0}^{{k\choose 2}}{{k\choose 2}\choose r}2^{{k\choose 2}-r}(-1)^r.
\end{eqnarray*}
The first factor inside the last summation is the number of ways to choose a subgraph $K_k|_{E'}$ with $r$ edges, and the second factor is the number of ways to choose a subgraph $K_k|_{E}$ containing $K_k|_{E'}$. Thus, we obtain finally
$$C_k
=\frac1{k!}\sum_{r=0}^{{k\choose 2}}{{k\choose 2}\choose r}(-2)^{-r}=\frac1{k!}\Bigl(1-\frac12\Bigr)^{{k\choose 2}}=\frac1{k!2^{{k\choose 2}}}.$$
This is what we claimed.

\subsection{Specialization to Abel polynomial}

We have proved already that the Abel polynomial possesses the binomial property and is
therefore umbral. Hence it suffices to find the sequence $\{a_n\}$
defined as
$$
a_n=\sum_G \frac{[q_n] A_G}{|\Aut~ G|},
$$
where the summation is carried over all connected graphs with~$n$ vertices.
As above, the last summation can be replaced by summation over all spanning subgraphs
of the complete graph~$K_n$. The contribution of any spanning tree into this
sum is $2^{(n-1)(n-2)/2}\cdot n$, since there are~$n$ ways to choose a root in the tree,
and each tree enters exactly $2^{(n-1)(n-2)/2}$ spanning (and necessarily connected) subgraphs of~$K_n$
(this is the number of ways to add other edges to the tree to obtain a spanning subgraph
containing this tree). Since, by Caley's theorem, the number of spanning trees in~$K_n$
is $n^{n-2}$, we obtain the desired equality
$$
a_n=2^{{n-1\choose2}}n^{n-1}.
$$

\end{document}